\documentclass{article}

\def\a{\alpha}

\def\fb{{\mathbf{f}}}

\newtheorem{thm}{Theorem}[section]
\newtheorem{lem}[thm]{Lemma}
\newtheorem{cor}[thm]{Corollary}

\newtheorem{rem}{Remark}[section]
\newtheorem{ex}{Example}[section]

\begin{document}

\title{On infinite dimensional Volterra type operators}

\author{Farrukh Mukhamedov$^1$, Mansoor Saburov$^2$}



\maketitle

\footnotetext[1]{Department of Comput. \& Theor. Sci., Faculty of
Sciences, IIUM, P.O. Box, 141, 25710, Kuantan, Pahang, Malaysia,
{\it e-mail:} {\tt far75m@yandex.ru} }

\footnotetext[2]{ Department of Mechanics \& Mathematics, National
University of Uzbekistan, Vuzgorodok, 100174 Tashkent, Uzbekistan }

 In this paper we study Volterra type
operators on infinite dimensional simplex. It is provided a
sufficient condition for Volterra type operators to be bijective.
Furthermore it is proved that the condition is not necessary. \vskip
0.3cm \noindent {\it
Mathematics Subject Classification}: 15A51, 47H60, 46T05, 92B99.


\section{Introduction}

Since Lotka and Volterra's seminal and pioneering
works (see \cite{V1}) many decades ago, modeling of interacting,
competing species have received considerable attention in the
fields of biology, ecology, mathematics (for example, see \cite{HS,T}). In their remarkably simple
deterministic model, Lotka and Volterra considered two coupled
nonlinear differential equations that mimic the temporal evolution
of a two-species system of competing predator and prey
populations. They demonstrated that coexistence of both species
was not only possible but inevitable in their model. Moreover,
similar to observations in real populations, both predator and
prey densities in this deterministic system display regular
oscillations in time, with both the amplitude and the period
determined by the prescribed initial conditions.
Note that in \cite{B,Ga,K,Lyu} finite dimensional Volterra and more general quadratic operators
were studied.

When a system is large enough, it is interesting to investigate quadratic Volterra operators define on an infinite dimensional simplex. First studies in this direction were considered in \cite{MAT}. Iterations of such operators define more complicated nonlinear operators. To better understanding the dynamics of such operators, it is important to study such nonlinear operators. The aim of this paper is to study more general class of nonlinear operators which contains a particular case of that quadratic Volterra operators. It is provided a sufficient
condition for Volterra type operators to be bijective

\section{Volterra type operators}
In this section we give definition of Volterra type operators and
study its basic properties. We prove that such operators are
bijective under certain condition.

Let
$$\ell_1=\{x=(x_n):\|x\|_1=\sum\limits_{k=1}^{\infty}|x_k|<
\infty\}$$ be the space of all absolutely convergent sequences.
The set $$S=\{x=(x_n)\in \ell_1: x_i\ge 0,
\sum\limits_{n=1}^{\infty}x_n=1 \}$$ is called {\it an infinite
dimensional simplex}.

It is known \cite{Roy} that $S=\overline{convh(Extr(S))}$ form a
convex closed set, here
$$Extr(S)=\{e^{(n)}=(\underbrace{0,0,\cdots,0,1,}_n0,\cdots)\}$$
is the extremal points of the $S$ and $conv(A)$ is the convex hull
of a set $A.$

Let $\alpha\subset N$ be an arbitrary set. The set
$$S_{\alpha}=\{x\in S: x_k=0, \forall k\in N \backslash \alpha\}$$
is called a {\it face} of the simplex. It is clear that $S_\alpha
$ is also will be the simplex.

{\it A relatively interior}  of a face $S_\alpha $ is defined by
$$
riS_{\alpha}=\{x\in S_\alpha: x_k>0, \forall k\in \alpha\}.$$

In what follows we are interested in the following operator
$V:S\to \ell_1$ defined by
\begin{equation}\label{v1}
(Vx)_k=x_k(1+f_k(x)), \  \  k\in N, \ x\in S,
\end{equation}
where a mapping $\fb: x\in S\to
(f_1(x),f_2(x),\cdots,f_n(x),\cdots)\in \ell_1$.

One can see that
$$\left|\sum\limits_{k=1}^{\infty}x_kf_k(x)\right|\le\sum\limits_{k=1}^{\infty}|x_kf_k(x)|\le
\sum\limits_{k=1}^{\infty}|f_k(x)|<\infty$$ for any $x\in S$, that
means that the operator $V$ is well defined.

\begin{thm}\label{VT1} Let an operator $V$ be given defined by (\ref{v1}). The following
conditions are equivalent:
\begin{itemize}
    \item[(i)] The operator $V$ is continuous in $\ell_1$ and $V(S)\subset
    S$. Moreover, $V(riS_\a)\subset riS_\a$ for every
    $\a\subset N$.
    \item[(ii)] The mapping $\fb$ satisfies the following
    conditions:
    \begin{itemize}
        \item[$1^0$] $\fb$ is continuous on $\ell_1$ topology.
        \item[$2^0$] For every  $k\in N$ one holds $f_k(x)\ge -1$ for all $x\in S.$
        \item[$3^0$] For any $x\in S$ one has
        $\sum\limits_{k=1}^{\infty}x_kf_k(x)=0.$
        \item[$4^0$] For any $\alpha\subset N$ one has $f_k(x)>-1$ for all $x\in riS_\alpha$ and
        $k\in\alpha.$
    \end{itemize}
    \end{itemize}
\end{thm}

{\bf Proof} (i)$\Rightarrow$(ii). The continuity of $V$  implies
$1^0$. Take $x\in S$, then $V(x)\in S$ yields that
\begin{itemize}
\item[(a)] $(V(x))_k\geq 0$, \item[(b)]
$\sum_{k=1}^\infty(V(x)_k)=1$.
\end{itemize}

Hence, from (a) it follows that $x_k(1+f_k(x))\geq 0$, which
implies $2^0$. From (b) one has $$\sum_{k=1}^\infty
x_k+\sum_{k=1}^\infty x_kf_k(x)=1,$$ which immediately yields
$3^0$.

Let $x\in riS_\alpha$, then $V(x)\in riS_a$, which with (\ref{v1})
and $x_k>0$ for all $k\in\a$ implies that $f_k(x)>-1$ for all
$k\in\alpha$.

The implication (ii)$\Rightarrow$(i) is evident.

We say that an operator $V:S\rightarrow S$ defined by (\ref{v1}) is
{\it Volterra type} if one of the conditions of Theorem \ref{VT1} is
satisfied. The corresponding mapping $\fb$ is called {\it
generating} for $V$.  By $\mathcal{V}$ we denote the set of all
Volterra type operators.

From Theorem \ref{VT1} we have

\begin{cor} The set $\mathcal{V}$ is convex, and for any
$V_1,V_2\in\mathcal{V}$ one has $V_1\circ V_2\in\mathcal{V}$, here
$\circ$ means composition of operators.
\end{cor}

Recall that the point $x\in S$ is called {\it fixed point} of $V$,
if $Vx=x$ The set of all fixed points of $V$ is denoted by
$Fix(V).$ From Theorem \ref{VT1} we immediately get the following

\begin{cor}\label{v3} For any  Volterra type operator $V$ one holds
\begin{itemize}
\item[(i)] $Extr(S)\subset Fix(V)$; \item[(ii)] Restriction of $V$
to any face of the simplex is again Volterra type operator.
\end{itemize}
\end{cor}

Let us consider more particular case of a mapping $\fb$. Namely
one has the following

\begin{thm}\label{VT2}   Let $f_k:\ell_1\rightarrow R$ ($k\in N$) be linear
functionals. Then a mapping $\fb:S\to\ell_1$ defined by
$\fb(x)=(f_1(x),f_2(x),\cdots,f_n(x),\cdots)$ satisfies the
conditions $1^0-4^0$ iff one has
\begin{equation}\label{v2}
f_k(x)=\sum\limits_{i=1}^{\infty}a_{ki}x_i, \ \ \ \ k\in N
\end{equation} with
$$a_{ki}=-a_{ik}, \ \ |a_{ki}|\le1 \ \ \ for \ \ every \ \
k,i\in N$$
\end{thm}

Note that Volterra type operators with generating mappings given by
(\ref{v2}) are called {\it quadratic Volterra operators}. In
\cite{MAT} such quadratic Volterra operators have been studied, and
it was shown that any such kind of operator is a bijection of $S$.
In the case under consideration, i.e. for Voterra type operators we
could not state an analogous result.

\begin{ex} Let us consider 1-dimensional simplex, i.e.
$S^1=[0,1]$. Define a Volterra type operator $V:[0,1]\to [0,1]$ by
$$
V(x)=x(1-\sin\pi x), \ \ \ x\in[0,1].
$$ A direct inspection shows that the defined operator is not
bijection.
\end{ex}

Now we are interested in finding some sufficient conditions for
Volterra type operators to be bijective.

\begin{thm}\label{VB} Let $V$ be a Volterra type operator given by (\ref{v1}).
Let
\begin{equation}\label{5-0}
\sum\limits_{k=1}^\infty x_kf_k(y)+\sum\limits_{k=1}^\infty
y_kf_k(x)\le 0  \ \ \textrm{for every} \ \ x,y\in S.
\end{equation}
be satisfied. Then $V$ is a bijection of $S$
\end{thm}

{\bf Proof}. First let us  prove that the $V$ is an injection.
Suppose that there are two distinct elements $x,y$ such that
\begin{equation}\label{xy}
V(x)=V(y).
\end{equation}

Without loss of generality we may assume that $x_i>0,$ $y_i>0$,
$\forall i\in N.$ If it is not true, then there is a face
$S_\alpha$, for some subset $\alpha\subset N$ of $S$ such that
$x,y\in riS_\alpha$ that is $x_i>0,$ $y_i>0,$ $\forall i\in\alpha$.
According to Theorem \ref{VT1} we have $V(S_\alpha)\subset S_\alpha$
therefore, due to Corollary \ref{v3} we can restrict $V$ to
$S_\alpha$ which is Volterra type too.

Now from (\ref{xy}) one gets
$$x_k(1+f_k(x))=y_k(1+f_k(y)), \ \ \forall k\in N$$
or
$$x_k-y_k+x_kf_k(x)-y_kf_k(x)=y_kf_k(y)-y_kf_k(x), \ \ \forall k\in N$$

$$(x_k-y_k)(1+f_k(x))=-y_k(f_k(x)-f_k(y)), \ \ \forall k\in N.$$
Since $1+f_k(x)>0,$ $\forall x\in riS$ and $y_k>0,$ $\forall k\in
N$, then

\begin{equation}\label{xy1}
sgn(x_k-y_k)=-sgn(f_k(x)-f_k(y)), \ \ \forall k\in N.
\end{equation}
Hence
\begin{equation}\label{xy2}
(x_k-y_k)(f_k(x)-f_k(y))\le 0, \ \ \forall k\in N
\end{equation}
whence
\begin{equation}\label{xy3}
\sum\limits_{k=1}^\infty(x_k-y_k)(f_k(x)-f_k(y))\le 0.
\end{equation}
Note that the last series absolutely converges, since
\begin{eqnarray*}
\left|\sum\limits_{k=1}^\infty(x_k-y_k)(f_k(x)-f_k(y))\right|&\le&\sum\limits_{k=1}^\infty
|x_k-y_k||f_k(x)-f_k(y)|\\
&\le&
\sum\limits_{k=1}^\infty(x_k+y_k)\left(|f_k(x)|+|f_k(y)|\right)\\
&\le
&2\left(\sum\limits_{k=1}^\infty|f_k(x)|+\sum\limits_{k=1}^\infty|f_k(y)|\right)<\infty
\end{eqnarray*}
From (\ref{xy3}) one finds
\begin{equation}\label{xy4}
\sum\limits_{k=1}^\infty x_kf_k(x)-\sum\limits_{k=1}^\infty
x_kf_k(y)-\sum\limits_{k=1}^\infty y_kf_k(x)+
\sum\limits_{k=1}^\infty y_kf_k(y)\le 0.
\end{equation}

By means of condition $3^0$, (\ref{xy4}) can be rewritten by
\begin{equation}\label{xy5}
\sum\limits_{k=1}^\infty x_kf_k(y)+\sum\limits_{k=1}^\infty
y_kf_k(x)\ge 0.
\end{equation}
From (\ref{5-0}) and (\ref{xy5}) we conclude that (\ref{xy2}) is
true if and only if
$$
(x_k-y_k)(f_k(x)-f_k(y))= 0, \ \ \forall k\in N.$$ The equality
(\ref{xy1}) with the last equality implies that $x=y.$ Thus,
$V:S\rightarrow S$ is injective.

Now let us show that $V$ is onto.
Denote
\begin{eqnarray*}
&&A_1=\{[1,n]\subset {N}: n\in {N}\}, \\
&&A_2=\{a\subset [1,n]:\left|[1,n]\backslash a\right|\ge 2, n\in
{N}\},\\
&& A_3=\{b\subset {N}: a\subset b, \ a\in A_1\cup A_2, \
\left|{N}\backslash b\right|<\infty\},\\
&&A=A_1\cup A_2\cup A_3. \end{eqnarray*}

In $A$ define an order by inclusion, i.e. $a\le b$ means that
$a\subset b$ for $a,b\in A.$ It is evident that $A$ is a
completely ordering set. We will prove surjection of $V$ by using
transfer induction method with respect to the set $A.$ Clearly
that operator $V$ is a surjection on $S_{\{1\}}.$ Suppose that for
an element $a\in A$ the operator $V$ is a surjective on $S_b$ for
every $b<a.$ Let us show that the $V$ is a surjection on $S_a$ as
well. Assume that $V(S_a)\neq S_a.$ For the boundary $\partial
S_a$ of $S_a$ we have $\partial S_a=\bigcup\limits_{c\in A: \
c<a}S_c.$ According to the assumption of the induction we get
\begin{equation}\label{xy6}
V(\partial S_a)=\partial S_a.
\end{equation}

On the other hand, there exist $x,y\in riS_a$ such that $x\in
V(S_a), \ y\notin V(S_a).$ The segment $[x,y]$ contains at least one
boundary point $z$ of the set $V(S_a).$ Since $V:S_a\rightarrow
V(S_a)$ is continuous and bijection, then the boundary point goes to
boundary one. Therefore, for $z\in riS_a, V^{-1}(z)\in \partial
S_a,$ which contradicts to (\ref{xy6}).

\begin{rem} Note that  the functionals (\ref{v2}) described
in Theorem \ref{VT2} satisfy the condition (\ref{5-0}).
\end{rem}

\section{Cubic Volterra operators}

From the previous section one arises a question: is there
non-trivial (except for linear functionals $f_k$) examples of
Volterra type operators for which condition (\ref{5-0}) is
satisfied. In this section  we introduce so called cubic Volterra
operators, and provide an examples of such operators which satisfy
condition (\ref{5-0}). Moreover, we establish that that condition is
indeed sufficient, i.e. an example of bijective cubic Volterra
operator will be given, which does not satisfy (\ref{5-0}).

Recall that a mapping $V:S\rightarrow S$ is called {\it cubic
stochastic operator} (shortly {\it c.s.o.}) if it is defined by
\begin{equation}\label{Cub}
(Vx)_k=\sum\limits_{i,j,l=1}^\infty p_{ijl,k}x_ix_jx_l, \ \ \forall
k\in N
\end{equation}
for all $x\in S$, here
\begin{equation}\label{Cubp}
\sum\limits_{k=1}^\infty p_{ijl,k}=1, \ \ \textrm{and} \ \
p_{ijl,k}=p_{\pi(i)\pi(j)\pi(l)}\ge 0,
\end{equation}
where $\pi$ is any permutation of the index set $\{i,j,l\}.$

Note that cubic stochastic operators were studied in
\cite{HS},\cite{RH}.

We  say that a c.s.o. $V$ is called {\it cubic Volterra operator
(c.V.o.)} if any face of the simplex is invariant with respect to
$V.$ One can prove the following

\begin{thm}\label{Cub1}
A c.s.o. $V$ is cubic Volterra operator if and only if
$p_{ijl,k}=0$ whenever $k\notin \{i,j,l\}.$ Moreover, any cubic
Volterra operator is Volterra type, and it can be represented by
\begin{equation}\label{Cub2}
(Vx)_{k}=x_k\left(x_k^2+3x_k\sum\limits_{i\in {N}_k}p_{ikk,k}x_i+
3\sum\limits_{i\in {N}_k}p_{iik,k}x_i^2+ 6\sum\limits_{i,j\in
{N}_{k},i<j}p_{ijk,k}x_ix_j\right),
\end{equation} for all $k\in N$, where ${N}_k={N}\backslash
\{k\}$.
\end{thm}

From this Theorem we immediately get the following

\begin{cor} If for a c.V.o. $V$ the coefficients $\{p_{ijl,k}\}$
do not depend one of indexes $\{i,j,l\}$, then $V$ becomes a
quadratic Volterra operator.
\end{cor}

Now we are going to provide an example of c.V.o. which satisfies
(\ref{5-0}).

{\bf Example 3.1.} Let us consider c.V.o. $V_0:S\rightarrow S$
defined by
$$p_{ikk,k}=1, \  p_{iik,k}=0, \  p_{ijk,k}=\frac13$$
here $i,j,k$ are pairwise distinct. Then $V$ has the following
form
\begin{equation}\label{ex1}
(V_0x)_k=x_k\left(x_k^2+3x_k\sum\limits_{i\in {N}_k}x_i+
2\sum\limits_{i,j\in {N}_k, i<j}x_ix_j\right),
\end{equation}
here, as before, ${N}_k={N}\backslash \{k\}$.

Now let us check (\ref{5-0}). One can see that
\begin{eqnarray*}
f_k(x)=x_k-\sum\limits_{i=1}^\infty x_i^2.
\end{eqnarray*}
Therefore, one gets
\begin{eqnarray*}
\sum\limits_{k=1}^\infty y_kf_k(x)+\sum\limits_{k=1}^\infty
x_kf_k(y)
&=&2\sum\limits_{k=1}^\infty x_ky_k-\sum\limits_{i=1}^\infty
x_i^2-\sum\limits_{i=1}^\infty y_i^2\leq 0,
\end{eqnarray*}
which implies (\ref{5-0}).

We are going to show that the condition (\ref{5-0}) is sufficient.
Next we give an example of bijective  c.V.o. which doesn't satisfy
that condition.

{\bf Example 3.2.} Let us consider an operator $V:S\rightarrow S$
defined by
\begin{equation}\label{ex2}
\left\{
\begin{array}{lll}
(Vx)_1=x_1^3\\[3mm]
(Vx)_2=x_2(x_2^2+3x_1)\\[3mm]
(Vx)_3=x_3(x_3^2+3(x_1+x_2)-3x_1x_2)\\[3mm]
\cdots\cdots\cdots\cdots\cdots\cdots\cdots\cdots\cdots\\
(Vx)_k=x_k\left(x_k^2+3\sum\limits_{i=1}^{k-1}x_i-3\sum\limits_{i,j=1,\
i<j}^{k-1}x_ix_j\right)\\
\cdots\cdots\cdots\cdots\cdots\cdots\cdots\cdots\cdots \\[3mm]
\end{array}
\right.
\end{equation}

Let us show that $V$ maps the simplex into itself.

The relations
\begin{equation}\label{ex22}
\sum\limits_{i=1}^nx_i-\sum\limits_{i,j=1,\
i<j}^nx_ix_j=\sum_{k=1}^{n-1}x_k
\bigg(1-\sum\limits_{i=k+1}^nx_i\bigg)+x_n\ge0 \end{equation} imply
that $(Vx)_k\ge0$ for all $k\in N$.

\begin{lem}\label{ex3} Let
$$W_k(x)=\left(\sum\limits_{i=k}^\infty x_i\right)^3+
3\sum_{i=1}^{k-1}x_i\left(\sum\limits_{i=k}^\infty x_i\right)^2+
3\sum\limits_{i,j=1,\ i\le j}^{k-1}x_ix_j\sum\limits_{i=k}^\infty
x_i, \ \ k\in N.$$ Then for all $k\in N$ one has
$$W_k(x)=(Vx)_k+W_{k+1}(x).$$
\end{lem}

{\bf Proof}
\begin{eqnarray*}
W_k(x)&=&\left(x_k+\sum\limits_{i=k+1}^\infty x_i\right)^3+
3\sum_{i=1}^{k-1}x_i\left(x_k+\sum\limits_{i=k+1}^\infty
x_i\right)^2\\
&&+ 3\sum\limits_{i,j=1,\ i\le
j}^{k-1}x_ix_j\left(x_k+\sum\limits_{i=k+1}^\infty x_i\right)\\
&=&x_k^3+3x_k^2\sum_{i=1}^{k-1}x_i+3x_k\sum_{i=1}^{k-1}x_i\sum\limits_{i=k+1}^\infty
x_i+3x_k\sum\limits_{i,j=1,\ i\le j}^{k-1}x_ix_j\\
&&+\left(\sum\limits_{i=k+1}^\infty x_i\right)^3+
3\left(x_k+\sum_{i=1}^{k-1}x_i\right)\left(\sum\limits_{i=k+1}^\infty
x_i\right)^2\\
&&+3\left(\sum\limits_{i,j=1,\ i\le
j}^{k-1}x_ix_j+x_k\sum_{i=1}^{k-1}x_i+x_k^2\right)\sum\limits_{i=k+1}^\infty
x_i\\
&=&x_k\left(x_k^2+3\sum_{i=1}^{k-1}x_i\sum\limits_{i=k}^\infty
x_i+3\sum\limits_{i,j=1,\ i\le j}^{k-1}x_ix_j\right)+W_{k+1}(x)\\
&=&(Vx)_k+W_{k+1}(x).
\end{eqnarray*}

From Lemma \ref{ex3} we obtain
\begin{eqnarray*}
1=\left(\sum\limits_{i=1}^\infty x_i\right)^3=(Vx)_1+W_2(x)
=(Vx)_1+(Vx)_2+W_3(x)=\cdots=\sum\limits_{i=1}^\infty (Vx)_i
\end{eqnarray*}

Now we show that the operator (\ref{ex2}) is bijective. To this end,
first we prove that the operator is injective. Indeed, let $x,y\in
S$ and $x\neq y.$ Then there exists $k_0\in N$ such that
$$x_{k_0}\neq y_{k_0}, \ \ \  x_i=y_i , \ \forall i=\overline{1,k_0-1}.$$

From (\ref{ex2}) we find that
$$
(Vx)_i=(Vy)_i, \ \ \forall
i=\overline{1,k_0-1}.
$$

From (\ref{ex22}) and $x_i=y_i $ $\forall i=\overline{1,k_0-1}$ it
follows that
$$
C:=\sum\limits_{i=1}^{k_0-1}x_i-\sum\limits_{i,j=1,\
i<j}^{k_0-1}x_ix_j=
\sum\limits_{i=1}^{k_0-1}y_i-\sum\limits_{i,j=1,\
i<j}^{k_0-1}y_iy_j\ge 0.$$

Consider a function $g(t)=t^3+3C\cdot t$, which is strictly
increasing on the segment $[0,1].$ Therefore,  for $x_{k_0}\neq
y_{k_0}$ we get
$$
(Vx)_{k_0}=g(x_{k_0})\neq g(y_{k_0})=(Vy)_{k_0},$$ which means
$Vx\neq Vy.$

Since the restriction of the operator (\ref{ex2}) to any face of $S$
is again of type (\ref{ex2}), then by similar argument used in the
proof of Theorem \ref{VB}, one can establish that the operator
(\ref{ex2}) is surjective.

Let us show that (\ref{5-0}) is not satisfied. From (\ref{ex2}) one
sees that
$$f_k(x)=x_k^2+3\sum\limits_{i=1}^{k-1}x_i-3\sum\limits_{i,j=1,\
i<j}^{k-1}x_ix_j-1.$$ The inequality (\ref{5-0}) is not satisfied if
we put $x=e^{(1)}$,$y=e^{(2)}$, indeed
$$
\sum\limits_{k=1}^\infty
e^{(1)}_kf_k(e^{(2)})+\sum\limits_{k=1}^\infty
e^{(2)}_kf_k(e^{(1)})=-1+2=1>0.
$$

\section*{Acknowledgement} The F.M. thanks Research Endowment Grant B of
 International Islamic University Malaysia, SAGA Fund P77c by MOSTI
 through the Academy of Sciences Malaysia(ASM) and TUBITAK.

\end{document}